\documentstyle{amsppt}

\NoBlackBoxes

\document

\topmatter
\title       On  Grothendieck Subspaces\endtitle
\author      S.~S.~Kutateladze\endauthor
\date        February 21, 2005\enddate
\address     Novosibirsk\endaddress

\endtopmatter

\noindent
Let $X$ be a~Riesz space, and let $Y$ be a~Kantorovich
(or Dedekind complete Riesz) space $Y$ with
base a~ complete Boolean algebra~$B$.
In~[1] we have described every order bounded operator $T$
from $X$ to $Y$ that may be presented as the difference of some Riesz homomorphisms.
An operator $T$ is such a~ difference if and only if
the kernel of the stratum $bT$ of~$T$
is a Riesz subspace of~$X$ for all $b\in B$.
We will give an analogous
description for an order bounded operator
$T$ whose modulus may be presented as the sum of two
Riesz homomorphisms in terms of the properties of the strata of~$T$.

Note that the sums of Riesz homomorphisms were first described
by  S.~J.~Bernau, C.~B.~Huijsmans, and B.~de~Pagter in terms of
$n$-disjoint operators in~[2].
A survey of some conceptually close results on $n$-disjoint operators
is given in~[3, \S 5.6].  In this note we reveal the connection between
the $2$-disjoint operators and Grothendieck subspaces.

\proclaim{Theorem}
The modulus of an order bounded  operator $T: X \to Y$
is the sum of some pair of Riesz homomorphisms
if and only if the kernel of each stratum
 $bT$ of~$T$ with $b\in B$ is a Grothendieck subspace of
 the ambient Riesz space~$X$.
\endproclaim

Recall that a subspace $H$ of a Riesz space is a~
{\it $G$-space\/} or {\it Grothendieck subspace\/}
provided that $H$ enjoys the following property:
$$
(\forall x,y\in H)\ (x\vee y\vee 0 + x\wedge y\wedge 0 \in H).
\tag1
$$

The history of~ (1) is as follows: In 1955 A.~Grothendieck distinguished
the subspaces  that satisfy (1) in the space $C(Q,\Bbb R)$
of continuous real functions on a compact space $Q$
(cp.~[4]). He determined
such a~subspace as 
the set of functions $f$ satisfying some family
$\roman A$ of relations of the form
$$
f(q_{\alpha}^1)=\lambda_\alpha f(q_{\alpha}^2)\
(q_{\alpha}^1,q_{\alpha}^2\in Q;\
\lambda_{\alpha}\in\Bbb R,\
\alpha\in{\roman A}).
$$
These subspaces due to A.~Grothendieck gave examples of
the $L_1$-predual Banach spaces other than
$AM$-spaces.
In 1969  J.~Lindenstrauss and D.~Wulpert characterized the
subspaces by using (1) and introduced the
concept of~$G$-space (cp.~[5]). Some related properties
of Grothendieck spaces are  presented  also in
[6] and~[7].

As in~[1], we rest the proof of the above-stated theorem
on the technique of ~``nonstandard scalarization''
reducing operator problems to the case of functionals.
Using the facts of Riesz space theory and Boolean
valued analysis we will presume the terminology and notation
of~[8].

We argue further as follows:
Using the functors of canonical embedding and descent to the
Boolean valued universe $\Bbb V^{(B)}$,
we reduce the matter to characterizing
a Grothendieck hyperspace that serves as
the kernel of an~order bounded functional
over a~dense subring of the reals~$\Bbb R$.
The scalar case  is settled by the following four lemmas.

\proclaim{Lemma~1}
A~functional $l$ is the sum of some pair of
Riesz homomorphisms if and only if
$l$ is a~ positive functional with kernel a~
Grothendieck subspace.
\endproclaim

\proclaim{Lemma~2}
Given an order bounded functional $l$ on a vector lattice~$X$,
assume that $l_+\neq 0$ and~$l_-\neq 0$.
The kernel $\ker(l)$ is a~Grothendieck subspace of if and only if
$l_+$ and~$l_-$ are Riesz homomorphisms on~$X$
(or, which is the same, $\ker(l)$ is a Riesz subspace).
\endproclaim

\proclaim{Lemma~3}
The kernel of an order bounded functional is
a~Grothendieck subspace if and only if so is the kernel
of the modulus of this functional.
\endproclaim

\proclaim{Lemma~4}
The kernel of an order bounded functional is
a~Grothendieck subspace if and only if the modulus of this
functional is the sum of a pair of~ Riesz homomorphisms.
\endproclaim

\enddocument